\newif\ifBASIC
\BASICfalse
\newif\ifWP
\WPfalse

\WPtrue

\documentclass{article}
\usepackage{amsmath,amsfonts,amssymb,latexsym,graphicx}

\ifWP

\makeatletter

\newif\iftwodates
\twodatesfalse

\renewcommand\maketitle{\begin{titlepage}%
  \let\footnotesize\small
  \let\footnoterule\relax
  \let \footnote \thanks
  \null\vfil
  \vskip 30\p@
  \begin{center}%
    {\LARGE \bf \@title \par}%
    \vskip 3em%
    {\large
     \lineskip .75em%
     \begin{tabular}[t]{c}%
       \@author
     \end{tabular}\par}%
     \vskip 1.5em%
  \end{center}\par
  \vfill
  \begin{center}
    \raisebox{1.5cm}{\includegraphics[width=0.58\textwidth]%
      {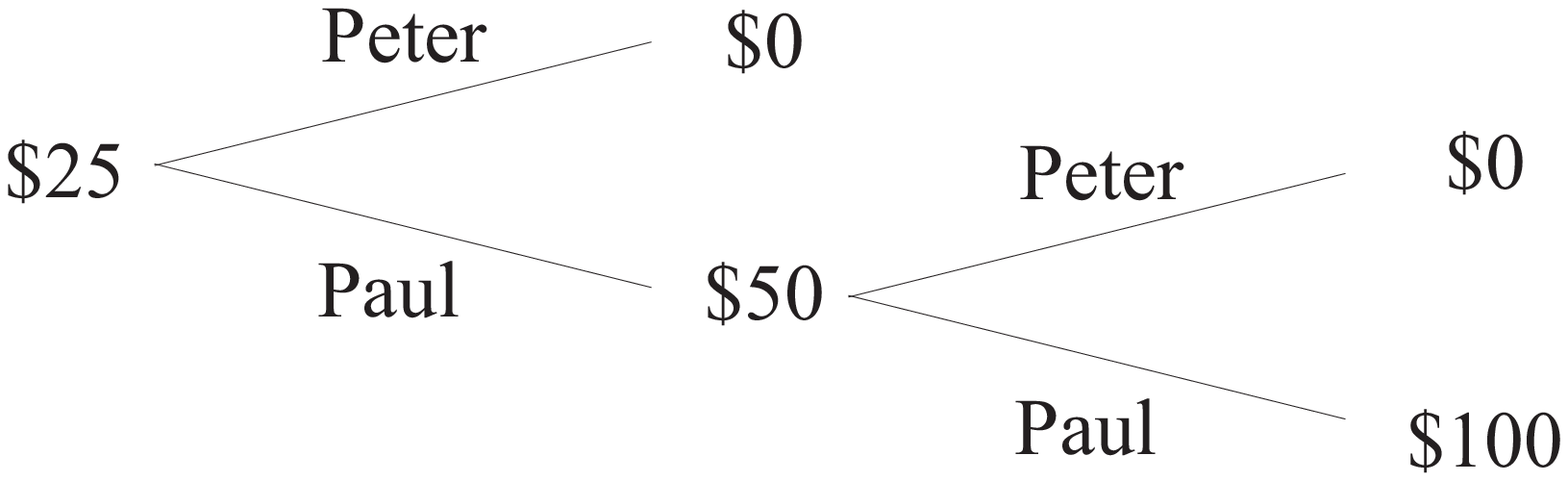}}%
    \hskip 3em%
    \includegraphics[width=0.29\textwidth]%
      {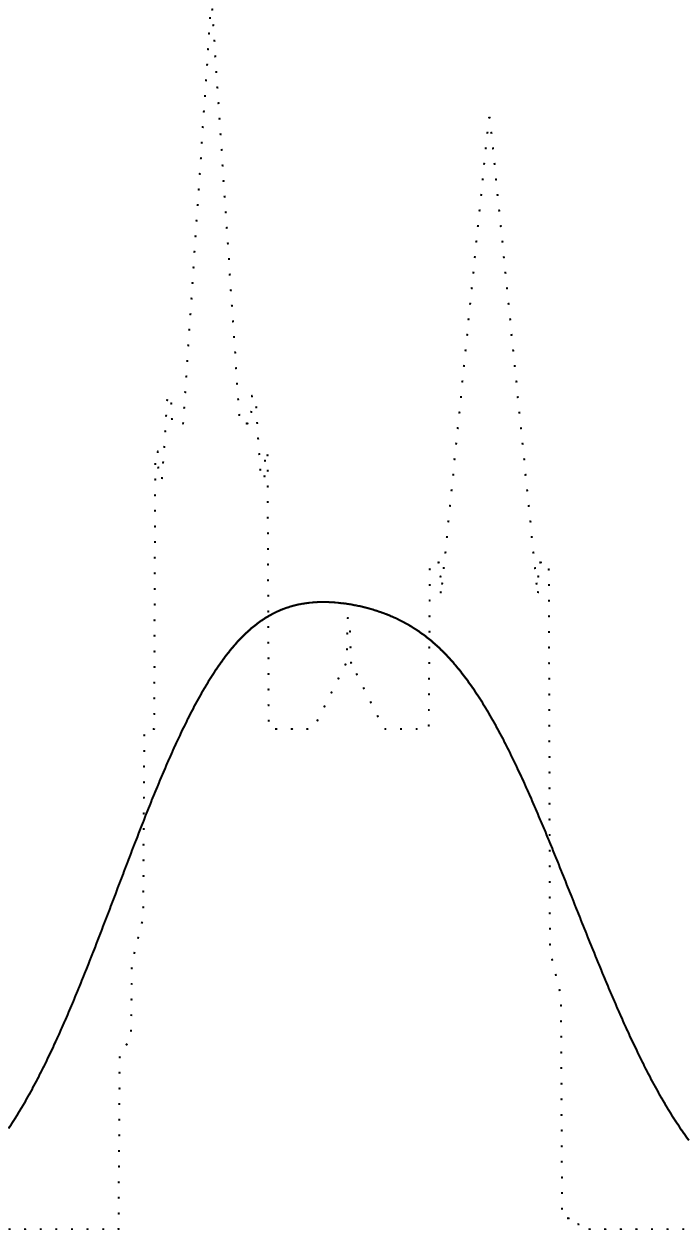}%
  \end{center}
  \@thanks
  \vfill
  \begin{center}
    {\large \bf The Game-Theoretic Probability and Finance Project}
  \end{center}
  \begin{center}
    {\large Working Paper \#\No}
  \end{center}
  \begin{center}
    {\iftwodates\large First posted \firstposted.
    Last revised \@date.\else\large\@date\fi}
  \end{center}
  \begin{center}
    Project web site:\\
    http://www.probabilityandfinance.com
  \end{center}
  \end{titlepage}%
  \setcounter{footnote}{0}%
  \global\let\thanks\relax
  \global\let\maketitle\relax
  \global\let\@thanks\@empty
  \global\let\@author\@empty
  \global\let\@date\@empty
  \global\let\@title\@empty
  \global\let\title\relax
  \global\let\author\relax
  \global\let\date\relax
  \global\let\and\relax
}

\renewenvironment{abstract}{%
  \titlepage
  \null\vfil
  \@beginparpenalty\@lowpenalty
  \begin{center}%
    \Large \bfseries \abstractname
    \@endparpenalty\@M
  \end{center}}%
  {\par\vfill\tableofcontents\endtitlepage}

\makeatother

\fi

\newcommand{\Prob}{\mathrm P}
\newcommand{\Expect}{\mathrm E}

\newcommand{\Bel}{\mathrm{Bel}}

\newcommand{\Q}{\mathrm Q}

\newcommand{\bbbr}{\mathbb{R}}		
\newcommand{\K}{\mathcal{K}}		

\newcommand{\SSS}{\mathcal{S}}

\newlength{\IndentI}
\newlength{\IndentII}
\newlength{\IndentIII}
\newlength{\WidthI}
\newlength{\WidthII}
\newlength{\WidthIII}
\ifBASIC
  \title{A Betting Interpretation \\
      for Probabilities and \\
      Dempster-Shafer Degrees of Belief\\}
  \author{Glenn Shafer%
  \thanks{Professor in the 
  Rutgers Business School - Newark and New Brunswick,
  180 University Avenue, Newark, New Jersey 07102 USA,
  and in the Department of Computer Science,
  Royal Holloway, University of London, Egham, Surrey TW20 0EX UK.
  E-mail gshafer@rutgers.edu;
  web site www.glennshafer.com.}}
\fi

\ifWP
  \title{A Betting Interpretation for Probabilities and
     Dempster-Shafer Degrees of Belief}
  \author{Glenn Shafer}
  \newcommand{\No}{31}
\fi

\begin{document}

\maketitle

\begin{abstract}
There are at least two ways to interpret numerical degrees of belief in terms of 
betting:
\begin{enumerate} 
   \item
      You can offer to bet at the odds defined by the degrees of belief.
   \item
      You can make the judgement that a strategy for taking advantage of such 
      betting offers will not multiply the capital it risks by a large factor.
\end{enumerate}
Both interpretations can be applied to ordinary additive probabilities and used to justify updating by conditioning. 
Only the second can be applied to Dempster-Shafer degrees of belief and used to justify Dempster's rule of combination.
\end{abstract}

\ifBASIC
  \vfill

  \pagebreak

  \tableofcontents

  \pagebreak
\fi 

\section{Introduction}

The meaning of numerical probability has long been a matter of contention.  
Sim\'eon Denis Poisson (1781--1840) distinguished between
objective and subjective probabilities \cite{poisson:1837}.  One recent 
philosophical introduction to probability lists five competing interpretations:
classical, frequency, propensity, logical, and subjective \cite{galavotti:2005}.

The classical and subjective interpretations both involve betting.
In the classical interpretation, the probability of an event is the correct price 
for a payoff that will equal one monetary unit if the event happens and zero
otherwise.  In the subjective interpretation, it is the price an individual
is willing to pay for this payoff.

This article explains another betting interpretation of probability.  Here I call it the \textit{Ville interpretation}, in recognition of Jean Andr\'e Ville (1910--1989), who first formulated it in his book on collectives \cite{ville:1939}.
Probabilities are prices under the Ville interpretation, 
just as they are under the classical and subjective interpretations.
But instead of asserting that these prices are correct in some unspecified sense
(as in the classical interpretation) or that some individual will pay them
(as in the subjective interpretation), we assert that no strategy for taking advantage of 
them will multiply the capital it risks by a large factor.  
The Ville interpretation derives from an older interpretation of probability, 
neglected in the English-language literature, which I call the 
\textit{Cournot interpretation} after Antoine Augustin Cournot (1801--1877).  According to the Cournot interpretation, the meaning of a probabilistic theory lies in the predictions that it makes with high probability.

As I explain in this article, the Ville interpretation can be applied both to ordinary additive probabilities and to the non-additive degrees of belief of the Dempster-Shafer calculus of belief functions.  It works for Dempster-Shafer degrees of belief in ways that the subjective interpretation does not.

\section{The Ville interpretation}

This section reviews how the Ville interpretation emerges
from older ideas and how it extends probability theory beyond its classical domain to games where the probabilities given and prices offered fall short of defining a probability distribution for all events of interest.
In Section~\ref{subsec:cournot}, I review briefly the history of the 
Cournot interpretation of ordinary 
probabilities.  In Section~\ref{subsec:ville}, I explain how the Ville interpretation  
is related to the Cournot interpretation.  In Section~\ref{subsec:forecast}, I illustrate
the power of the Ville interpretation using the example of probability forecasting,
and in Section~\ref{subsec:game}, I explain its role more generally in game-theoretic
probability.

\subsection{Cournot}\label{subsec:cournot}

The standard procedure for testing a probabilistic theory involves picking out an event to which
the theory gives very small probability:  we reject the theory if the event happens.
In fact, this seems to be the only way to test a probabilistic theory.
Because Cournot was the first to state that mathematical probability
makes contact with phenomena only by ruling out events given very small probability          
(\cite{cournot:1843}, p.~58), 
the prediction that 
\begin{equation}\label{eq:cournotclass}
    \text{an event of very small probability will not happen}
\end{equation}
has been called \textit{Cournot's principle}.
In the first half of the twentieth century, many European scholars, including \'Emile Borel, Paul L\'evy, Maurice Fr\'echet, and Andrei Kolmogorov, contended that Cournot's principle is fundamental to the meaning and use of mathematical probability
\cite{shafer/vovk:2006}.  As Borel said, we evoke ``the only law of chance'' when we single out an event of very small probability and predict it will not happen.  (Or when, equivalently, we single out an event of very high probability and predict that it will happen.)  Let us call the thesis that such predictions constitute the meaning of probability the \textit{Cournot interpretation of probability}.  

Cournot, Fr\'echet, and Kolmogorov are often called frequentists.  This is misleading.  These authors did believe that 
the probability of an event will be approximated by the frequency with which it happens in independent trials, but they considered this ``law of large numbers'' a consequence of Cournot's principle together with Bernoulli's theorem, which gives very high probability to the approximation holding.  The true frequentists, such as John Venn, saw no sense in Bernoulli's theorem; probability is frequency, they believed, and so it is silly to try to prove that frequency will approximate probability \cite{venn:1888}.  

Of course, events of very small probability do happen.  An experiment may have a very large number of possible outcomes, each of which has very small probability, and one of which must happen.  So Cournot's principle makes sense only if we are talking about particular events of very small probability that are salient for some reason:  perhaps
because they are so simple, perhaps because they have high probability under a plausible
alternative hypothesis, or perhaps simply because they were specified in advance.
There may be a substantial number of events that are salient in this way, but this is not 
a problem if we set our threshold for small probability low enough, because the 
disjunction of a number of events with very small probably will still have reasonably small probability.

In order to put the Cournot interpretation into practice, we must also decide how small a probability we can neglect.  This evidently depends on the context.  Borel distinguished between what was negligible at the human level, at the terrestrial level, and at the cosmic level \cite{borel:1939}.

In using the Cournot interpretation, we must also bear in mind its role in testing and giving meaning to a probabilistic theory as a whole.  Strictly speaking, it gives direct meaning only to probabilities that are very small (the event will not happen) or very large (the event will happen).  It gives no meaning to a probability of $40\%$, say.  But when a probabilistic theory says that many successive events are independent and all have probability $40\%$, it gives probabilities close to one for many aspects of this sequence of events.  Probabilistic theories in which probabilities evolve (stochastic processes) also give probabilities close to one to many statements concerning what happens over time, so they can also be tested and acquire meaning by Cournot's principle.

Although it was widely accepted in continental Europe in the middle of the twentieth century, the Cournot interpretation never gained a significant foothold in the English-language literature, and awareness of it receded as English became the language of science and mathematics after World War II.  We find only isolated affirmations of it after about 1970.  In the article on probability in the Soviet Mathematical Encyclopedia, for example, we find the assertion that only probabilities close to zero or one have empirical meaning \cite{prokh/sevast:1987}.  For more on the history of the Cournot interpretation see \cite{martin:1996,martin:2003,martin:2005,shafer:2007,shafer/vovk:2006}.

\subsection{From Cournot to Ville}\label{subsec:ville}

When a probability distribution is used to set betting odds, there is a well known
relationship between the happening of events of small probability and the success of 
betting strategies.  
The event that a given betting strategy multiplies the capital it risks by $1/\alpha$ or more 
has probability $\alpha$ or less.
Conversely, for every event of probability $\alpha$ or less there is a bet
that multiplies the capital it risks by $1/\alpha$ or more if the event happens.  
So it is natural to consider, as an alternative to 
Cournot's principle, the principle that
\begin{equation}\label{eq:cournotgame}
    \text{a strategy will not multiply the capital it risks by a large factor.}
\end{equation}
Let us call this \textit{Ville's principle}.  Let us call the thesis that predictions of the form~(\ref{eq:cournotgame}) constitute the meaning of probability the \textit{Ville interpretation} of probability.

Ville's principle is equivalent to Cournot's principle whenever a probability distribution is given for the events being considered and the two principles are made specific, with the specific event and small probability mentioned in Cournot's principle matching the specific strategy and large factor mentioned in Ville's principle.  But when the two principles are considered more abstractly, without $\alpha$ and the particular event or strategy being specified, they differ in two important respects:
\begin{enumerate}
\item
Ville's principle gives us more guidance than Cournot's principle.  It tells us to specify a strategy for betting, not merely a single event of small probability. We found it necessary to elaborate Cournot's principle by saying that the event of very small probability should be specified in advance.  The corresponding coda for Ville's principle is also needed, but it is less easily overlooked, because a betting strategy cannot be implemented unless it is specified in advance.
\item
Ville's principle has a broader scope than Cournot's principle.  Cournot's principle applies only when there is a probability distribution for the events under discussion. 
Ville's principle applies whenever prices for gambles are given, even if these prices fall short of defining probabilities for events.
\end{enumerate}

To see some of the implications of Ville's principle giving us more guidance, consider how testing is usually implemented.  A test of a probabilistic theory usually begins with a test statistic, say $T(y)$, where $y$ is an outcome that is 
to be observed.  If the theory specifies a probability distribution $\Prob$ for $y$, then we reject the theory at the significance level $\alpha$ when we observe a value $y$ such that
$$
                     T(y) \ge c,
$$ 
where $c$ is a number such that $\Prob \{ T(y) \ge c \} \le \alpha$.  Ville's principle tells us to implement this idea in 
a particular way:  our test statistic is the capital $\K(y)$ achieved by a specified betting strategy that starts with some initial capital $\K_0$ and does not risk losing more than $\K_0$.  We reject the theory at the significance level $\alpha$
when we observe a value $y$ such that
$$
         \K(y) \ge \K_0/\alpha.
$$
Markov's inequality tells us that $\Prob \{ \K(y) \ge \K_0/\alpha \} \le \alpha$.%
\footnote{In general, Markov's inequality says that a nonnegative random
variable $X$ satisfies 
  $$
     \Prob \left( X \geq  \Expect(X) / \alpha \right) \leq \alpha.
  $$
Because the betting strategy uses the odds set by $\Prob$, the expected value of the final capital $\K(y)$ is the initial capital $\K_0$.  Because the strategy risks only the initial capital, the final capital $\K(y)$ cannot be negative.}

When we adopt a betting strategy with which to test a probability distribution $\Prob$, we are implicitly specifying an alternative hypothesis $\Q$ that we can plausibly adopt if we reject $\Prob$.  To see that this is so, let us suppose, for simplicity, that $\K_0=1$ (the strategy risks one unit of capital), and that there are only finitely or countably many possible values for $y$.  In this case, we can define $\Q$ by
\begin{equation}\label{eq:likrat}
         \Q(y) := \K(y) \Prob(y).
\end{equation}
It is easy to see that $\Q$ is a probability distribution:  (1) 
$\Q(y)\ge 0$ because $\Prob(y)$ is a probability
and $\K(y)$ is the final capital for a betting strategy that does not 
risk its capital becoming negative, and (2) $\sum_y \Q(y)=1$ because it is the expected
payoff under $\Prob$ of a gamble that costs one unit.
Equation~(\ref{eq:likrat}) tells us that the final capital 
$\K(y)$ is the likelihood ratio $\Q(y)/\Prob(y)$, a measure of how much
the observed outcome $y$ favors $\Q$ over $\Prob$.

\subsection{Probability forecasting}\label{subsec:forecast}

As a first example of how Ville's principle and the Ville interpretation apply even when prices offered fall short of defining a probability distribution $\Prob$ for all events of interest, consider a game in which a forecaster announces probabilities successively, observing the outcome of each preceding event before giving the next probability:

\bigskip

\noindent
\textsc{Probability Forecasting Game}

\parshape=6
\IndentI  \WidthI
\IndentI  \WidthI
\IndentII \WidthII
\IndentII \WidthII
\IndentII \WidthII
\IndentII \WidthII
\noindent
$\K_0:=1$.\\
FOR $n=1,2,\dots,N$:\\
  Forecaster announces $p_n\in [0,1]$.\\
  Skeptic announces $s_n \in \bbbr$.\\
  Reality announces $y_n \in \{0,1\}$.\\
  \hbox to \WidthII
  {$\K_n := \K_{n-1} + s_n (y_n - p_n)$.
  \refstepcounter{equation}\label{eq:capital}
  \hfil (\theequation)}

\bigskip

\noindent
This is a perfect-information game; the three players move in sequence, and they 
all see each move as it is made.  The game continues for $N$ rounds.

The number $p_n$ can be thought of as the price of a ticket that pays 
the amount $y_n$.  
Skeptic can buy any number $s_n$ of the tickets.  Since he pays $p_n$ for each ticket
and receives $y_n$ in return, his net payoff is $s_n(y_n - p_n)$.
The number $s_n$ can be positive or negative.  By choosing $s_n$ positive,
Skeptic buys tickets; by choosing $s_n$ negative, he sells tickets.

Within the game, the $p_n$ are simply prices.  But we think of them as
Forecaster's probabilities:  $p_n$ is Forecaster's probability that Reality
will choose $y_n=1$.  On the other hand, Forecaster need not have a  
joint probability distribution $\Prob$ for Reality's moves $y_1,\dots,y_N$.
He simply chooses $p_n$ as he pleases at each step.

Skeptic tests Forecaster's $p_n$ by trying to increase
his capital using them as prices.  If Skeptic succeeds---i.e., if he makes
$K_N$ large without risking more than his initial capital $\K_0$, 
then we conclude that Forecaster is not
a good probability forecaster.  Ville's principle says that if Forecaster is a good
forecaster, then Skeptic will not achieve a large value for his final capital $K_N$ 
without risking more than $\K_0$.

What does it mean for Skeptic not to risk more than $\K_0$?  It means that his moves 
do not allow Reality to make his final capital $\K_N$ negative.  Since Reality
can always keep Skeptic from making money (by choosing $y_n=0$
if $s_n$ is positive and $y_n=1$ if $s_n$ is negative), she can make $\K_N$ negative
as soon as Skeptic lets $\K_n$ become negative for any $n$.  So in order to deny
Reality the option of making $\K_N$ negative, Skeptic must choose each $s_n$ so as to
deny Reality the option of making $\K_n$ negative.  By~(\ref{eq:capital}), this means 
choosing $s_n$ in the interval
\begin{equation}\label{eq:safe}
   -\K_{n-1}/p_n \le s_n \le \K_{n-1}/(1-p_n).
\end{equation}
For brevity, let us say that Skeptic plays \textit{safely} if he always chooses
$s_n$ satisfying (\ref{eq:safe}), and lets us call a strategy for Skeptic \textit{safe}
if it always prescribes $s_n$ satisfying~(\ref{eq:safe}).

We can get back to classical probability by assuming that Forecaster follows 
a strategy based on a joint probability distribution $\Prob$
for $y_1,\dots,y_N$ and perhaps other events outside the game, the strategy 
being to set $p_n$ equal to $\Prob$'s conditional probability for $y_n=1$
given what has been observed so far.  But Ville's principle is powerful even in the absence of a specified strategy
for Forecaster.  It is all we need in order to derive
various relations, such as the law of the large numbers, the law of the iterated
logarithm, and the central limit theorem, that classical probability theory 
says will hold 
between the probabilities $p_1,\dots,p_N$ and the outcomes $y_1,\dots,y_N$.  
It turns out, for example, that Skeptic can play safely in such a way 
that either the relative frequency 
of $1$s among $y_1,\dots,y_N$,
$
          \sum_{n=1}^N y_n/N,
$
approximates the average probability forecast,
$
          \sum_{n=1}^N p_n/N,
$
or else $\K_N$ becomes very large
(\cite{shafer/vovk:2001}, p.~125).  Because it tells us that 
$\K_N$ will not become very large very large, Ville's principle therefore
implies that 
$
          \sum_{n=1}^N y_n/N
$
will approximate
$
          \sum_{n=1}^N p_n/N.
$
This is a version of the law of large numbers.

\subsection{Game-theoretic probability}\label{subsec:game}

Probability forecasting is only one example where prices fall short
of defining a probability distribution.  In many other examples, the shortfall is 
substantially greater.

One class of such examples arises in finance theory, where the price for a security
at the beginning of the day can be thought of as the price for a ticket that pays
what the security is worth at the end of the day.  Here the roles of Forecaster and 
Reality are both played by the market that sets the prices, and the role of Skeptic
is played by a speculator.  Over a period of $N$ days, they play a perfect-information
game much like our Probability Forecasting Game:

\bigskip

\noindent
\textsc{Market Game}

\parshape=6
\IndentI  \WidthI
\IndentI  \WidthI
\IndentII \WidthII
\IndentII \WidthII
\IndentII \WidthII
\IndentII \WidthII
\noindent
$\K_0:=1$.\\
FOR $n=1,2,\dots,N$:\\
  Market announces opening price $p_n\in [0,\infty)$.\\
  Speculator announces $s_n \in \bbbr$.\\
  Market announces closing price $y_n \in [0,\infty)$.\\
  $\K_n := \K_{n-1} + s_n (y_n - p_n) $.

\bigskip
\noindent
Here $s_n > 0$ when Speculator goes long in the security, and
$s_n < 0$ when he goes short. 

A cornerstone of finance theory is the \textit{efficient market hypothesis}, which
states that a speculator cannot expect to make money using publicly available 
information.  Efforts to formulate this hypothesis more precisely usually 
start with the questionable assumption that market prices are governed, in some
sense, by a probability distribution.  Ville's principle offers an 
alternative way of making the hypothesis precise:  we can say that Speculator 
will not make $\K_N$ large while playing safely.  This version of the efficient 
market hypothesis can be tested directly, without making any probabilistic 
assumptions \cite{GTP23}.  It also
implies a number of stylized facts about financial markets, including 
the $\sqrt{dt}$ effect \cite{GTP5} and the relation between the volatility
and average of simple returns called the CAPM \cite{vovk/shafer:2008}.

Shafer and Vovk \cite{shafer/vovk:2001} give other examples of games where 
prices fall short of defining a probability distribution.  It turns out that 
many of the usual results of probability theory can be extended to such games,
provided that we adopt Ville's principle.  In general, we call 
the study of such games \textit{game-theoretic probability}.

The results in \cite{shafer/vovk:2001} are concerned with strategies for Skeptic or Speculator in a probability game; they say that this player can multiply their capital by a large factor if some result in probability theory or finance theory does not hold.  It is also fruitful, however, to consider how Forecaster or Market can play against such strategies for Skeptic or Speculator.  It turns out that they can do this effectively, and this gives a new method of making predictions, called \textit{defensive forecasting} \cite{vovk/takemura/shafer:2005, vovk/shafer:2005}.

\section{The judgement of irrelevance in updating by conditioning}\label{sec:conditioning}

How should Forecaster's probabilities change when he learns new information?  

An important school of thought, called \textit{Bayesian} in recent decades, contends that when we learn $A$, we should update our probability for $B$ from $\Prob(B)$ to 
\begin{equation}\label{eq:updating}
              \frac{\Prob(A \&B)}{\Prob(A)}.
\end{equation}
The change is called \textit{conditioning}.  Bayesians acknowledge that it is appropriate only if we judge $A$ to be the only relevant information we have learned (\cite{definetti:1970},
Section~11.2.2, \cite{bernardo/smith:1994}, p.~45).%
\footnote{The authors just cited, de Finetti and Bernardo and Smith, go on to say that irrelevance usually fails; when we learn $A$ we usually learn other information that will also modify our judgement concerning $B$.  Nevertheless, updating by~(\ref{eq:updating}) is widely taught and implemented.} 

In this section, I review arguments for the updating rule~(\ref{eq:updating}), with attention to how they account for the judgement of relevance and irrelevance.  I consider the argument originally given by Abraham De Moivre, the variation given by Bruno de Finetti, and another variation that is based on Ville's principle.  Only the argument from Ville's principle uses the judgement of relevance.

\subsection{De Moivre's argument}\label{subsec:moivrecond}

Abraham De Moivre was the first to state the \textit{rule of compound probability}.  In the second edition of his \textit{Doctrine of Chances}, published in 1738 \cite{demoivre:1718}, he stated the rule as follows:
\begin{quotation}
\noindent
  \dots the Probability of the happening of two Events dependent, is the product of the Probability
  of the happening of one of them, by the Probability which the other will have of happening, when the
  first shall have been consider'd as having happen'd\dots
\end{quotation}
This rule can be written
\begin{equation}\label{eq:rulecp}
   \Prob(A \& B) = \Prob(A) \Prob(B|A),
\end{equation}
where $\Prob(A \& B)$ is the probability of the happening of $A$ and $B$, $\Prob(A)$ is the probability of the happening of $A$, and $\Prob(B|A)$ is the probability which $B$ will have of happening, when $A$ shall have been consider'd as having happen'd.  

The twentieth century abandoned De Moivre's way of talking about probabilities.  Now we call $\Prob(B|A)$ the \textit{conditional probability} of $B$ given $A$, and we say that it is defined by the equation 
\begin{equation}\label{eq:defcp}
   \Prob(B | A) := \frac{\Prob(A \& B)}{\Prob(A)},
\end{equation}
provided that $\Prob(A)\neq 0$.  This makes~(\ref{eq:rulecp}) a trivial consequence of a definition.  But for De Moivre, (\ref{eq:rulecp}) was more substantive.  It was a consequence of how probability is related to price.

De Moivre gave an argument for the rule of compound probability on pp.~5--6 of his second edition.  He used a language that is somewhat unfamiliar today; he talked about the values of gamblers' expectations.  But it is true to his thinking to say that the probability of an event is the price (or the fair price, if you prefer) for a ticket that pays $1$ if the event happens and $0$ if it does not happen.  (An expectation is the possession a ticket with a uncertain payoff, and its value is the price you should pay for the ticket.)  Using the language of tickets, payoffs, and price, we can express his argument as follows:
\begin{enumerate}
  \item
The price of a ticket that pays $1$ if $A$ happens is $\Prob(A)$.
  \item
Assume one can buy or sell any number of such tickets, even fractional amounts.  So $\Prob(A) x$ is the price of a ticket that pays $x$ if $A$ happens, where $x$ is any real number.  (Buying a negative amount means selling.)
  \item
After $A$ happens (or everyone learns that $A$ has happened and nothing else), $\Prob(B|A)$ is the price of a ticket that pays $1$ if $B$ happens.
  \item
So starting with $\Prob(A) \Prob(B|A)$, you can get $1$ if $A \& B$ happens.  You use the $\Prob(A) \Prob(B|A)$ to buy a ticket that pays $\Prob(B|A)$ if $A$ happens, and then, if $A$ does happen, you use the $\Prob(B|A)$ to buy a ticket that pays $1$ if $B$ also happens.
  \item
So $\Prob(A) \Prob(B|A)$ is the value of a ticket that pays $1$ if $A \& B$ happens.
\end{enumerate}     

De Moivre's argument is unconvincing to modern readers because we do not accept his starting point---his unexamined assumption that an expectation has a well defined numerical value.  Our positivist heritage demands that such numbers be cashed out in some way that can be observed.

\subsection{De Finetti's version of the argument}\label{sec:finetti}

Bruno de Finetti (1906--1985) had a way of responding to the positivist challenge.  For him, probability is specific to an individual.  An individual's probability for an event $A$ is the price the individual sets for a ticket that returns $1$ if $A$ happens---the price at which he is willing to trade in such tickets, buying or selling as the occasion arises.

As for the conditional probability $\Prob(B|A)$, de Finetti proposed a betting interpretation that avoids references to a situation after $A$ has happened or is known to have happened.  For him, $\Prob(B|A)$ is the price of a conditional ticket---the price of a ticket that pays $1$ if $B$ happens, with the understanding that the transaction is cancelled (the price is refunded and no payoff is made if $B$ happens) if $A$ does not happen.

With these interpretations, de Finetti was able to formulate a version of De Moivre's argument that leaves aside the notion of changing probabilities.  We situate ourselves at the beginning of the game, as it were, and argue as follows:
\begin{enumerate}
  \item
$\Prob(A)$ is the price at which I am willing to buy or sell tickets that pay $1$ if $A$ happens.
  \item
I am willing to buy or sell any number of such tickets, even fractional amounts.  So $\Prob(A) x$ is the price I will pay for a ticket that pays $x$ if $A$ happens, where $x$ is any real number.
  \item
$\Prob(B|A)$ is the price I am willing to pay for a ticket that pays $1$ if $B$ happens, with the understanding that this price is refunded if $A$ does not happen.
  \item
It follows that I am willing to pay $\Prob(A) \Prob(B|A)$ to get back $1$ if $A$ and $B$ both happen.  You can prove this by selling me two tickets:
\begin{itemize}
\item
  For $\Prob(A) \Prob(B|A)$, a ticket that pays $\Prob(B|A)$ if $A$ happens.
\item
  For $\Prob(B|A)$, a ticket that pays $1$ if $B$ and $A$ both happen, with the price being refunded if $A$ does not happen.
\end{itemize}
If $A$ and $B$ both happen, I end up with $1$, less the $\Prob(A) \Prob(B|A)$ I paid for the first ticket; the payoff from the first ticket is cancelled by the cost of the second.  If $A$ does not happen, I lose only the $\Prob(A) \Prob(B|A)$, the second purchase having been cancelled.  If $A$ happens but $B$ does not, I again lose only the $\Prob(A) \Prob(B|A)$, the cost of the second purchase being cancelled by the payoff on the first.
  \item
So $\Prob(A) \Prob(B|A)$ is the price I am willing to pay for $1$ if $A \& B$ happens---i.e., my probability for $A \& B$.
\end{enumerate}     
As a coda, we may add de Finetti's argument for the price being unique.  De Moivre had taken it for granted that the value of a thing is unique.  De Finetti, using his assumption that we are willing to buy and sell any amount, argued that we must make the probability unique in order to prevent an opponent from extracting an indefinite amount of money from us.

De Finetti's version of the argument comes closer to modern mathematical rigor than De Moivre's, because it leaves aside the notion of something being ``consider'd as having happen'd'', for which De Moivre gave no set-theoretic exegesis.  But some such notion must still be used in order to extend the argument to a justification for using conditional probabilities as one's new probabilities after something new is learned.  We must explain why the price $\Prob(B|A)$ for the conditional ticket on $B$ given $A$ should not change when $A$ and nothing else is learned.  There is a large literature on how convincingly this argument can be made; some think it requires that a protocol for new information be fixed and known in advance.  See \cite{shafer:1985} and references therein.

\subsection{Making the argument from Ville's principle}

Ville's principle, like Cournot's, can usually be applied directly only to a run of events, in which a strategy has time to multiply the capital it risks substantially (or, in the case of Cournot's principle, we can identify an event of very small probability).  So in order to apply Ville's principle to the problem of changing probabilities that are neither very small nor very large, we must imagine them being embedded in a longer sequence of similar probabilities for similar events.  This is how probability judgments are often made:  we judge that an event is like an event in some repetitive process for which we know probabilities \cite{shafer/tversky:1985}.

In de Finetti's picture, we make a probability judgement $\Prob(A)=p$ by saying that $p$ is the price at which we are willing to buy or sell tickets that pay $1$ if $A$ happens.  (I omit needed caveats:  that we buy and sell only to people who have the same knowledge as ourselves, that this is only the price we might be inclined to set if we were inclined to gamble, etc.)  In Ville's picture, we make a probability judgement $\Prob(A)=p$ by saying that if we do offer such bets on $A$, and on a sequence of similar events in similar but independent circumstances, then an opponent \textit{will not succeed in multiplying the capital they risk in betting against us by a large factor}.  Let us abbreviate this to the statement that an opponent \textit{will not beat the probability}.

In this terminology, our task is to show that the following claim holds:
\begin{equation}\label{eq:claim}
\begin{minipage}{.8\linewidth}
\noindent
Suppose we are in a situation where we judge that an opponent will not beat $\Prob(A)$ and $\Prob(A\& B)$.  Suppose we then learn $A$ and nothing more.  Then we can include $\Prob(A\& B)/\Prob(A)$ as a new probability for $B$ among the probabilities that we judge an opponent will not beat.
\end{minipage}
\end{equation}
In one respect, we are following De Moivre more faithfully than de Finetti did.  De Finetti's mathematical argument is concerned only with prices in a single situation.  Here we propose, like De Moivre, to give an argument that relates prices over two situations:  an initial situation and a subsequent situation where our additional knowledge is $A$ and nothing more.  This is normal for the game-theoretic framework reviewed in Sections~\ref{subsec:forecast} and~\ref{subsec:game}; there we apply Ville's principle to games with many rounds. 

Here is the argument for~(\ref{eq:claim}) from Ville's principle:
\begin{enumerate}
  \item
An opponent will not beat the probabilities $\Prob(A)$ and $\Prob(A \& B)$.  This means that a strategy for the opponent that buys and sells tickets on $A$ and $A \& B$ at these prices, along with similar tickets on other events, will not multiply the capital risked by a large factor.
  \item
We need to show that this impossibility of multiplying the capital risked still holds for strategies that are also allowed to use $\Prob(A\& B)/\Prob(A)$ as a new probability for $B$ after $A$ and nothing more is known.
  \item
It suffices to show that if $\SSS$ is a strategy against all three probabilities ($\Prob(A)$ and $\Prob(A\& B)$ in the initial situation and $\Prob(A\& B)/\Prob(A)$ later), then there exists a strategy $\SSS'$ against the two probabilities ($\Prob(A)$ and $\Prob(A \& B)$ in the initial situation) alone that risks no more capital and has the same payoffs as $\SSS$.
  \item
Let $M$, which may be positive or negative or zero, be the amount of $B$ tickets $\SSS$ buys after learning $A$.  To construct $\SSS'$ from $\SSS$, we delete this purchase of $B$ tickets and add 
\begin{equation}\label{eq:addtickets}
   M \text{ tickets on } A \& B \quad 
       \text{and} \quad
   -M\frac{\Prob(A \& B)}{\Prob(A)} \text{ tickets on } A
\end{equation}
to $\SSS$'s purchases of tickets on $A$ and $A \& B$ in the initial situation.  
\begin{itemize}
\item
  The tickets in~(\ref{eq:addtickets}) have zero net cost: 
$$
  M \Prob(A \& B) - M\frac{\Prob(A \& B)}{\Prob(A)} \Prob(A) = 0.
$$
So $\SSS'$ uses the same capital in the initial situation as $\SSS$.
\item
  The payoffs  of the tickets in~(\ref{eq:addtickets}) are the same as the net payoffs of the $M$ tickets deleted from $\SSS$:
\begin{align*}
 & 0 \text{ if } A \text{ does not happen};\\
 & -M \frac{\Prob(A\&B)}{\Prob(A)} \text{ if } A \text{ happens but not } B;\\
 & M\left(1-\frac{\Prob(A\&B)}{\Prob(A)}\right) \text{ if }  A \text{ and } B \text{ both happen}.
\end{align*}
so $\SSS'$ uses no more capital than $\SSS$ after the initial situation and has the same payoffs in the end.
\end{itemize}
\item
By hypothesis, $\SSS'$ will not multiply the capital it risks by a large factor.  So $\SSS$, which risks the same capital and has the same payoffs, does not either.
\end{enumerate}
See \cite{shafer/gillett/scherl:2003} for an extension of this argument to Peter Walley's updating principle for upper and lower previsions.

\subsection{The judgement of irrelevance}

The argument from Ville's principle for using conditional probability as one's new probability uses the role and implications of knowledge in a way that de Finetti's argument does not.  
\begin{itemize}
\item
De Finetti argued for the conditional probability $\Prob(B|A)$ being the price in the initial situation for a conditional purchase---a purchase of a $B$ ticket on the condition that $A$ happens.  He then merely asserted, with no argument, that it should remain the price for this purchase after we learn that $A$ happens and nothing more.  
\item
The Ville argument, in contrast, is truly an argument for $\Prob(B|A)$ being the price for a $B$ ticket in the new situation where we have learned that $A$ happened and nothing more. 
\end{itemize}
The Ville argument is able to bring knowledge into the story because it looks what can be accomplished by different strategies.  What a strategy can accomplish depends on what information is available.

It is important to understand how the caveat ``nothing more'' enters into the Ville argument.  The argument depends on constructing a strategy $\SSS'$ for the initial situation alone that is equivalent to a strategy $\SSS$ that makes additional bets in the later situation where $A$ is known.  If something more than $A$ is known, and $\SSS'$ uses this additional information as well ($\SSS$'s purchase of the $M$ tickets depends on it), then the construction is not possible.

We can of course relax the requirement that nothing more be known than $A$'s happening.  The essential requirement is that nothing more be known that can help an opponent multiply his capital.  In this case, we may say that the happening of $A$ is our only \textit{relevant information}.  We may have learned many other things by the time or at the time when we learned $A$, but none of them can provide further help to a strategy for betting against the probabilities.

\section{Judgements of irrelevance in the Dempster-Shafer calculus}\label{sec:dempster}

The Dempster-Shafer theory of belief functions extends conditional probability to a calculus for combining probability judgements based on different bodies of evidence.  Judgements of irrelevance enter into this calculus explicitly and pervasively.  These judgements can be explained in terms of Ville's principle in the same way as the judgement of irrelevance in the case of updating by conditioning on $A$:  they are judgements that once certain information is taken into account, other information is of no help to a strategy for betting against certain probabilities.

In this section, I list the Ville judgements of irrelevance required by various operations in the Dempster-Shafer calculus (Section~\ref{subsec:oper}), and I discuss how attention to these judgements in applications can strengthen the calculus's usefulness (Section~\ref{subsec:appl}).

\subsection{Basic operations}\label{subsec:oper}

The Dempster-Shafer calculus derives from a series of articles by A.~P.~Dempster, recently republished along with other classic articles on the calculus in \cite{yager/liu:2008}.  The calculus was described in detail in \cite{shafer:1976} and reviewed in \cite{dempster:2008}.  Without reviewing the examples and details readers can find in these references, I give here an overview of four related operations:  the transfer of belief, conditioning, independent combination, and Dempster's rule of combination.  In each case, I explain the judgement of irrelevance involved.

I omit two other important operations, natural extension and marginalization, because they do not require judgements of irrelevance.

\paragraph*{Transfer of belief.}
Suppose $\mathbf{X}$ is a variable, whose possible values form the set $\mathcal{X}$, and suppose $\Prob$ is a probability distribution on $\mathcal{X}$, expressing our probability judgements about the value of $\mathbf{X}$.  Suppose $\omega$ is another variable, with the set of possible values $\Omega$, for which we do not have a probability distribution.  Suppose further that  $\Gamma$ is a multivalued mapping from $\mathcal{X}$ to $\Omega$ (a mapping from  $\mathcal{X}$ to non-empty subsets of $\Omega$).  Then we can define a function $\Bel$ on subsets of $\Omega$ by setting 
\begin{equation}\label{eq:beldef}
    \Bel(A) := \Prob\{x | \Gamma(x) \subseteq A \}.
\end{equation}
A function defined in this way is called a \textit{belief function}.  We call $\Bel(A)$ its \textit{degree of belief} in $A$.

We can give $\Bel$'s degrees of belief a Ville interpretation under the following conditions:
\begin{enumerate}
\item
The probability distribution $\Prob$ has a Ville interpretation:  no betting strategy will beat the probabilities it gives for $\mathbf{X}$.
\item
The multivalued mapping $\Gamma$ has this meaning:
\begin{equation}\label{eq:gamma}
     \text{If } \mathbf{X} = x, \text{ then } \omega\in\Gamma(x).
\end{equation}
\item
Learning the relationship~(\ref{eq:gamma}) between $\mathbf{X}$ and $\omega$ does not affect the impossibility of beating the probabilities for $\mathbf{X}$.  (This is the irrelevance judgement.)
\end{enumerate}
The Ville interpretation that follows from these conditions is one-sided:  a strategy that buys for $\Bel(A)$ tickets that pay $1$ if $\omega \in A$ (and makes similar bets on the strength of similar evidence) will not multiply the capital its risks by a large factor.

\paragraph*{Conditioning.}
Suppose we modify the preceding setup by allowing the subset $\Gamma(x)$ of $\Omega$ to be empty for some $x$.  In this case, condition~(\ref{eq:gamma}) tells us that the event $\{x | \Gamma(x) \neq \emptyset\}$ happened, and if we judge that we have learned nothing else that can help a strategy beat $\Prob$'s probabilities, then we are entitled to condition $\Prob$ on this event.  This results in replacing~(\ref{eq:beldef}) by 
$$
    \Bel(A) = \frac{\Prob\{x | \Gamma(x) \subseteq A \;\;  \&  \;\;
                               \Gamma(x) \neq \emptyset \}}{\Prob \{x | \Gamma(x) \neq \emptyset\}}.
$$
The judgements of irrelevance that justify this equation can be summarized by saying that aside from the impossibility of the $x$ for which $\Gamma(x)=\emptyset$, learning~(\ref{eq:gamma}) does not provide any other information that can help a strategy beat the probabilities for $\mathbf{X}$.

\paragraph*{Independent combination.}
Suppose $\Prob_1$ and $\Prob_2$ are probability distributions on $\mathcal{X}_1$ and $\mathcal{X}_2$, respectively, expressing our probability judgements about the values of the variables $\mathbf{X}_1$ and $\mathbf{X}_2$, respectively.  What judgement is involved when we say further that the product probability measure $\Prob_1 \times \Prob_2$ on $\mathcal{X}_1 \times \mathcal{X}_2$ expresses our probability judgement about $\mathbf{X}_1$ and $\mathbf{X}_2$ jointly?

This question is not answered simply by saying that $\mathbf{X}_1$ and $\mathbf{X}_2$ are probabilistically independent, because probabilistic independence, in modern probability theory, is a property of a joint probability distribution for two variables, not a judgement outside the mathematics that justifies adopting the product distribution as a joint probability distribution for them.

De Finetti's betting interpretation of probability does give an answer to the question:  we should adopt the product distribution if learning the value of one of the variables and nothing else will not change the prices we are willing to offer on the other variable.  

The Ville interpretation gives an analogous answer:  we should adopt the product distribution if we make the judgement that knowing the value of one of the variables and nothing more would not help a strategy beat the probabilities for the other variable.

Dempster-Shafer theory extends the idea of independent combination to belief functions, by considering two multivalued mappings, say a mapping $\Gamma_1$ from $\mathcal{X}_1$ to non-empty subsets of $\Omega_1$ and a mapping $\Gamma_2$ from $\mathcal{X}_2$ to non-empty subsets of $\Omega_2$.  Suppose $\Gamma_1$ and $\Gamma_2$ have these meanings, where $\omega_1$ and $\omega_2$ are variables that take values in $\Omega_1$ and $\Omega_2$, respectively:
\begin{equation}\label{eq:judge1}
     \text{If } \mathbf{X_1} = x, \text{ then } \omega_1\in\Gamma_1(x).
\end{equation}
\begin{equation}\label{eq:judge2}
     \text{If } \mathbf{X_2} = x, \text{ then } \omega_2\in\Gamma_2(x).
\end{equation}
Then we can form a belief function $\Bel$ for the pair $(\omega_1,\omega_2)$:
$$
    \Bel(A) = (\Prob_1 \times \Prob_2) \{(x_1,x_2) | \Gamma_1(x_1) \times \Gamma_2(x_2) \subseteq A\}
$$
for $A \subseteq \Omega_1 \times \Omega_2$.  To justify this, we must make the Ville judgement justifying the formation of the product distribution $ \Prob_1 \times \Prob_2 $ and also the judgement that learning~(\ref{eq:judge1}) and~(\ref{eq:judge2}) does not help beat the probabilities given by $\Prob_1$ or $\Prob_2$.  This goes beyond the individual judgements that learning~(\ref{eq:judge1}) does not help beat $\Prob_1$ and that learning~(\ref{eq:judge2}) does not help beat $\Prob_2$.

\paragraph*{Dempster's rule of combination.}
Dempster's rule concerns the combination of two bodies of evidence bearing on the same variable $\omega$.  Given the ideas we have just reviewed, it is most easily stated by considering two multivalued mappings from the different probability spaces to the same space $\Omega$, say $\Gamma_1$ from $(\mathcal{X}_1,\Prob_1)$ and $\Gamma_2$ from $(\mathcal{X}_2,\Prob_2)$.  They have the usual meaning:
\begin{equation}\label{eq:judge3}
     \text{If } \mathbf{X_1} = x, \text{ then } \omega\in\Gamma_1(x).
\end{equation}
\begin{equation}\label{eq:judge4}
     \text{If } \mathbf{X_2} = x, \text{ then } \omega\in\Gamma_2(x).
\end{equation}
Even if both $\Gamma_1(x_1)$ and $\Gamma_2(x_2)$ are always non-empty, their intersection may be empty.  When we learn~(\ref{eq:judge1}) and~(\ref{eq:judge2}), we learn that the event $\{(x_1,x_2) | \emptyset \neq \Gamma_1(x_1) \cap \Gamma_2(x_2)\}$ has happened.

Conditioning on the intersection being non-empty, we obtain the belief function $\Bel$ on $\Omega$ given by
$$
    \Bel(A) := \frac{(\Prob_1 \times \Prob_2 )\{(x_1,x_2) | \emptyset \neq \Gamma_1(x_1) \cap \Gamma_2(x_2) \subseteq A \}}
                   {(\Prob_1 \times \Prob_2 ) \{(x_1,x_2) | \emptyset \neq \Gamma_1(x_1) \cap \Gamma_2(x_2)\}}.
$$
In this case, the required Ville judgements are those involved in forming the product measure, along with the judgement that learning~(\ref{eq:judge1}) and~(\ref{eq:judge2}) does not help beat the probabilities given by the product measure aside from providing the information that $\{(x_1,x_2) | \emptyset \neq \Gamma_1(x_1) \cap \Gamma_2(x_2)\}$ has happened.

\subsection{Discussion}\label{subsec:appl}

In \cite{shafer:1976}, I stated that Dempster's rule of combination is appropriate when the bodies of evidence underlying individual belief functions are independent.  The Ville judgements I have just detailed elaborate this notion of independence, in a way that should be useful in applications.

In our various writings on belief functions and in debates with critics, A.~P.~Dempster and I frequently took the view that the notions of independence and conditioning involved in Dempster's rule are the same as in ordinary probability theory.  The analysis of this article vindicates this view in some degree, insofar as it has shown that the judgements of irrelevance required for Dempster's rule have the same general form as judgements of irrelevance that justify the formation of product measures in ordinary probability theory and updating by conditioning in Bayesian reasoning.  The analysis has also revealed, however, the complexity that can be involved in judgements of this general form.

The critics often demanded, of course, explanations of independence and conditioning that were consistent with de Finetti's explanation of the meaning of these concepts in the Bayesian calculus.  Here I have argued that de Finetti's explanations are not as convincing as sometimes thought even for Bayesian updating:  they justify the pricing of conditional tickets but not the changes in price from one state of knowledge to another.  In any case, they surely do not extend to the Dempster-Shafer case, where no embedding of the rules in a static picture seems to be possible.  For the process of combining evidence, we need a more dynamic picture, which is provided by the Ville interpretation.

It is easy to construct examples in which the Ville irrelevance judgements required for Dempster's rule are unreasonable or clearly wrong.  It is also easy enough to construct examples in which these judgements are reasonable; I gave some such examples in the 1980s (see for example \cite{shafer/tversky:1985}).  Existing applications of the Dempster-Shafer calculus would be enriched, however, by a systematic examination of the reasonableness of the irrelevance judgements they require.  A clearer understanding of these judgements might also help us construct Dempster-Shafer models for complex scientific problems where the irrelevance judgements need to justify ordinary probabilities and Bayesian reasoning seem unreasonably strong.

\bibliographystyle{plain}
\bibliography{betting}
\end{document}